\documentclass[11pt]{article}
\usepackage{amsmath,amssymb,amsthm,latexsym}
\newtheorem{thm}{Theorem}[section]

\newtheorem{lem}[thm]{Lemma}
\newtheorem{prop}[thm]{Proposition}
\newtheorem{prin}[thm]{Principle}

\newtheorem{remark}[thm]{Remark}
\newtheorem{conj}[thm]{Conjecture}
\newtheorem{ques}[thm]{Questions}

\newcommand{\thmref}[1]{Theorem~\ref{#1}}

\newcommand{\propref}[1]{Proposition~\ref{#1}}

\newcommand{\secref}[1]{Section~\ref{#1}}
\newcommand{\figref}[1]{Figure~\ref{#1}}
\newcommand{\bx}{\hfill$\Box$}

\renewcommand\a{\alpha}         

\newcommand\g{\gamma}

\newcommand\e{\epsilon}
\renewcommand\l{\lambda}
\renewcommand\L{\Lambda}
\newcommand\D{\Delta}
\newcommand\G{\Gamma}
\newcommand\f{\frac}
\newcommand\Z{\mathbb Z}
\newcommand\R{\mathbb R}
\newcommand\U{\mathbb H}
\newcommand\C{\mathbb C}
\newcommand\A{\mathbb A}
\newcommand\Q{\mathbb Q}
\renewcommand\Re{\mbox{Re~}}
\renewcommand\Im{\mbox{Im~}}
\newcommand\sech{\mbox{sech}}

\newcommand{\quo}[1]{SL_{{#1}}({\Z})\backslash SL_{{#1}}({\R})/
SO_{{#1}}({\R})}

\renewcommand\i{^{-1}}

\renewcommand\({\left(}         
\renewcommand\){\right)}

\include{epsf}

\begin{document}

\title{The Highest-Lowest Zero and other Applications of
Positivity}
\author{Stephen D. Miller\thanks{The author was supported by
 an NSF
Postdoctoral Fellowship during this work.} \\ {\scriptsize
Department of Mathematics , Yale University , P.O. Box 208283, New
Haven, CT 06520-8283} }
\date{May 16, 1999}
\vspace{-100 mm}
\maketitle

\vspace{-10 mm}
\begin{abstract}

The first nontrivial zeroes of the Riemann $\zeta$ function
are $\approx \f12\pm14.13472i$.  We investigate the question of whether or
not any other L-function has a higher lowest zero.  To do so we try to
quantify the notion that the L-function of a ``small'' automorphic
representation (i.e. one with small level and archimedean type) does not
have small zeroes, and vice-versa.  We prove that many types of
automorphic
L-functions have a lower first zero than $\zeta$'s (see
Theorems~\ref{realarch} and \ref{gl2thm}).  This is done using Weil's
explicit formula with carefully-chosen test functions.  When this method
does not immediately show L-functions of a certain type have low zeroes,
we then attempt to turn the tables and show no L-functions of that type
exist.  Thus the argument is a combination of proving low zeroes exist and
that certain cusp forms do not.  Consequently we are able to prove
vanishing theorems and improve upon existing bounds on the Laplace
spectrum on $L^2(\quo n)$.  These in turn can be used to show that
$\quo{68}$ has a discrete, non-constant, non-cuspidal eigenvalue outside
the range of the continuous spectrum on
$L^2(SL_{68}({\R})/SO_{68}({\R}))$, but that this never happens for $\quo
n$ in lower rank.  Another application is to cuspidal cohomology: we show
there are no cuspidal harmonic forms on $\quo n$ for $n < 27$.

\end{abstract}

\section{Introduction}

The Riemann $\zeta$ function's first critical zeroes are
 surprisingly large: about $ \f 12 \pm 14.13472i$.  Our main
interest in this paper is the following question:
\begin{quotation}
Does any other automorphic L-function have a larger first zero?
\end{quotation}
This question was raised by Odlyzko (\cite{Odlyzko}),
who proved that the Dedekind zeta function of any number field
has a zero whose imaginary part is less than 14.  Odlyzko also proved
related conditional results for Artin L-functions.

Every automorphic L-function conjecturally factors into products
of standard L-functions of cusp forms on $GL_n$ over the
rationals, and we shall be content to discuss
these.\footnote{Nevertheless our arguments work under various
wider assumptions, as they are mostly sensitive to the analytic
properties of the L-function.}  In fact, by twisting a cuspidal
automorphic representation of $GL_n/{\Q}$ by a power of the
determinant, it is possible to shift the zeroes any amount
vertically, so we restrict ourselves to studying cuspidal
automorphic representations $\pi_=\otimes_{p\le\infty}\pi_p$
 of $GL_n/{\Q}$ whose
central character is normalized to be trivial. In most examples
coming from number theory the archimedean type $\pi_\infty$ is
real, i.e. the gamma factors multiplying $L(s,\pi)$ have real
shifts. Our first result answers the question for such cusp forms:

\begin{thm}
\label{realarch}
 Let $\pi$
be a cuspidal automorphic representation of $GL_n$ over $\Q$
 with a real archimedean type and a trivial central
character.  Then
$L(s,\pi)$ has a low zero which either (i) is on the critical axis
 between $\f12 \pm 14.13472i$ or (ii) violates the generalized
Riemann
hypothesis (GRH) in an effective range.
\end{thm}

When we speak of a zero violating GRH ``in an effective
range,'' we mean that should conclusion (i) fail, then one could
theoretically find an effective constant $T>0$ such that the box
$(\f 12,1)\times [-T,T]i$ contains a zero.  For brevity we will use the
following terminology:

\vspace{.5cm}

{\bf Definition}: An L-function has a {\em low zero}
if it either vanishes on the critical axis between
 $\f12 \pm 14.13472i$, or violates GRH in an effectively-bounded
range (see \secref{wogrh}).

\vspace{.5cm}

We will use this definition to state unconditional results, but not much
is actually gained philosophically or numerically in this problem
by assuming GRH.

The L-functions in \thmref{realarch} include those of Dirichlet
characters, rational elliptic curves, and conjecturally all
rational abelian varieties. Of course they are also expected to
include all Artin L-functions, for example L-functions of Galois
representations. We have been unable to squeeze our technique to
answer Odlyzko's question in full generality, but can prove many
cases.  For example:

\begin{thm} Let $\pi$ be a cuspidal automorphic representation
of $GL_2$ over $\Q$ with a trivial central character.  Then $L(s,\pi)$
has a {\em low zero} (which is on the critical axis
 between $\f12 \pm 14.13472i$ or else violates GRH in an effective
range).
\label{gl2thm}
\end{thm}

This includes modular form and Maass form L-functions.

Other results can be proven about low zeroes.  For example, every
L-function which is related to itself by an odd functional
equation automatically vanishes at $s=1/2$. For a fixed degree
$n$, most cuspidal automorphic representations of $GL_n$ over $\Q$
with a trivial central character have low zeroes.  In fact, the
possible exceptions all lie in a bounded subset of the unitary
dual and have bounded level.  This subset tends to be devoid of
cusp forms, which is why our method is successful.  Thus Odlyzko's
question is related to vanishing theorems about automorphic forms.

Our technique uses Weil's explicit formula relating the coefficients
and zeroes of automorphic L-functions.  It is a variation
on the Stark-Odlyzko positivity technique, as formulated by
Serre, Poitou, Mestre, and others -- see \cite{Odlyzko} for a
survey.
In particular,
one can compute an exact formula for sum of certain test
functions over
the critical zeroes.  If we use a test function which is
positive only in a certain range,
then finding this sum is positive ensures a zero in that range.
On the other hand, if this sum
is negative, then we can often construct another test function
which
is positive in the critical strip, yet whose sum over the zeroes
is negative.  This contradiction shows that the L-function actually
could not have existed to begin with.
Our main difficulty is that it is often very difficult to construct this
second test function given the failure of the first.

The latter contradiction, of positive terms yielding a negative sum, can
be used to prove vanishing theorems about automorphic forms, since they
cannot exist when their L-functions do not.
Independent of our interest in low zeroes, this leads to
applications
in group cohomology and spectral theory.

\subsection*{Other applications}

One of the consequences of the Ramanujan-Selberg
 temperedness conjecture is that the discrete cuspidal spectrum
of the laplacian $\D$ on \newline $L^2(\quo n)$ is contained in
the continuous spectrum of $\D$ on $L^2(SL_n({\R})/SO_n({\R}))$.
(We always normalize $\D$ so that this continuous spectrum is the
interval $[\f{n^3-n}{24},\infty)$.) This consequence should be
true more generally for congruence covers of $\quo n$, but in this
particular case slightly more was proven in \cite{Miller}:

\begin{thm}(\cite{Miller}): There exists a constant $c>0$
such that the Laplace eigenvalue of every cusp form  $\phi$ on
$\quo n$ satisfies
$$\l(\phi)> \l_1(SL_n({\R})/SO_n({\R}))+cn.$$
\label{imrn}
\end{thm}

Our new result is superior for small $n$:
\begin{thm} Let $\phi$ be a cuspidal eigenfunction
of the non-euclidean
laplacian $\D$ on $SL_n({\Z})\backslash SL_n({\R})/SO_n({\R})$.
Then $\phi$'s Laplace eigenvalue satisfies
\begin{equation}
\l(\phi)>\f{n^3-4n}{24} + 25.92\(1+\f{1}{n-1}\).
\end{equation}
\label{newlamb1}
\end{thm}

It can be applied to answer a question of Alexander Lubotzky: when
does the eigenvalue of a noncuspidal, square-integrable
 eigenfunction of the laplacian
 on $SL_n({\Z})\backslash SL_n({\R})/SO_n({\R})$ lie outside
$[\f{n^3-n}{24},\infty)$?

\begin{thm} If~ $n\le 67$, any non-constant
eigenfunction of $\D$ in \newline $L^2(\quo n)$ has Laplace
eigenvalue greater than $\f{n^3-n}{24}$, but the first Laplace
eigenvalue of $SL_{68}({\Z})\backslash
SL_{68}({\R})/SO_{68}({\R})$ is in fact approximately
$$12906.6< \f{68^3-68}{24}=13098.5 .$$
\end{thm}
Finally, we can apply our technique to cuspidal cohomology and extend
a result in \cite{Miller}, where it was shown that
$SL_n({\Z})\backslash SL_n({\R})/SO_n({\R})$
has no {\em harmonic} cuspidal automorphic forms for $n<23$:

\begin{thm}
The constant-coefficients cuspidal cohomology of $SL_n({\Z})$
$$H_{cusp}^{\cdot}(SL_n({\Z});{\C})=0$$ vanishes for
$1<n<27.$
\end{thm}

The technique used to prove this theorem is related to the one
in \cite{Miller}.  Fermigier \cite{Fermig} had a similar, but weaker,
result using positivity with a different L-function.  Here
we combine both methods to go further.

{\bf Acknowledgements}: We wish to thank Don Blasius, William
Duke, Benedict Gross,
Alexander Lubotzky, Andrew Odlyzko, Ilya Piatetski-Shapiro, Vladimir
Rokhlin,
Peter Sarnak,
Jean-Pierre Serre,
Gunther Steil, Andrew Wiles, and Gregg Zuckerman for their discussions.
Our point of view on L-functions was influenced by the discussion
in \cite{Rud-Sar}.  Support was provided by National Science Foundation
Graduate and Postdoctoral Fellowships and a Yale Hellmann
fellowship during stays at Princeton University, Yale University, and
the University of California at San Diego.  All numerical computations
were made with Mathematica v.3 on an Intel Pentium II
300 MHz system running Windows NT 4.0 and Slackware Linux 2.0.30.

\section{L-functions and positivity}

 By conjectures of Langlands the most
general automorphic L-function is a product of
standard L-functions of cuspidal automorphic representations
$\pi=\otimes_{p\le\infty}\pi_p$
 on $GL_m$ over the rational adeles $\A_{\Q}$.  These
``primitive'' L-functions are
degree
$m$ Euler products
$$L(s,\pi)=\prod_{p\mbox{ prime}} \
\prod_{j=1}^m (1-\a_{p,j} p^{-s})^{-1}~,~
\a_{p,j}\in \C$$ and have completions
$$\L(s,\pi)=\prod_{j=1}^m \pi^{(-s+\eta_j)/2}
\Gamma\(\f{s+\eta_j}{2}\) L(s,\pi)~,~ \eta_j \in {\C}$$
which are entire unless $m=1$ and $L(s)=\zeta(s)$.  We have used
the duplication property of the gamma function in writing
the gamma factors in this way.  The conductor
is $D$, and for $\pi_p$ unramified, the $\a_{p,j}$ are
Hecke eigenvalue parameters and the $\eta_j$ are related to
the archimedean
parameters of $\pi_\infty$.
 With this normalization $\L(s,\pi)$ has the functional
equation
$$\L(s,\pi)=\tau_\pi D^{-s}\L(1-s,\tilde{\pi})~
,~\tau_\pi \in {\C}~,~|\tau_\pi|=\sqrt{D}~,~D>0,$$ where
$\tilde{\pi}$ is the contragredient representation to $\pi$. The
Jacquet-Shalika (\cite{JS}) bounds imply that
\begin{equation}
\Re{\eta_j}>-\f12.
\label{jsbound}
\end{equation}

\subsection{Weil's formula}

The explicit formula of Andr\'e Weil equates a sum over the zeroes of an
L-function with a sum over its coefficients and gamma factors:
\begin{equation}
\sum_{\L(\f{1}{2}+i\g,\pi)=0}h(\g) = 2 \Re\left\{ \sum_{j=1}^m
l(\eta_j) - \sum_{n=1}^{\infty}\f{c_n}{\sqrt{n}}g(\log n)\right\}+
g(0)\log{D},
\label{weil}
\end{equation}
where $g$ is an even, differentiable real function,
$$\hat{g}(r)=h(r)=\int_{{\R}}g(x)e^{irx}dx,$$
$$\Gamma_{{\R}}(s)=\pi^{-s/2} \Gamma(s/2),$$ and
$$l(\eta)=\f{1}{2\pi}\int_{{\R}}h(r)
\f{\Gamma_{{\R}}'}{\Gamma_{{\R}}}\(\f{1}{2}+\eta+ir\)dr$$
$$=\f{1}{2\pi}\int_{{\R}}h(r)\(\f{-\log\pi}{2}\)dr +
\f{1}{2\pi}\int_{{\R}}h(r)\f{\Gamma'}{2\G}\(\f{1}{4}+\f{\eta}{2}
+\f{ir}{2}\)dr$$
$$=-\f{\log\pi}{2}g(0) -
\f{1}{2}\int_0^{\infty}\(\f{g(x/2)
e^{-(1/4+\eta/2)x}}{1-e^{-x}}-\f{g(0)}{e^{x}x}\)dx.$$
Here we have made use of the fact that $L(s,\pi)$ is entire; for $\zeta(s)$
 and Rankin-Selberg
L-functions there is a polar term that will be introduced
when needed later on.  See \cite{Rud-Sar} for a proof of
(\ref{weil}).

If $g$ is supported in the interval $[-\log 2,\log 2]$ then the formula
can be viewed as giving the value of the sum over the zeroes from the
gamma
factors:
\begin{equation}
\sum h(\g) = 2 \Re \sum_{j=1}^m l(\eta_j) + g(0)\log D.
\label{simplexplic}
\end{equation}
The basis of the positivity technique is the observation that
 if $h(\g)
\ge 0$ for each zero, then the sum on the right-hand side
of (\ref{simplexplic}) must
 also be
positive.  This immediately gives a lower bound on the conductor
$D$, which is the original application of the
positivity technique.  Fortunately the sum
on the right-hand side of (\ref{simplexplic}) is explicitly computable in
 terms of the
$\eta_j$'s and $D$; if it is negative then the L-function
 $L(s,\pi)$
cannot exist and hence neither can the original cusp form $\pi$.

Upon assuming GRH, let
$$\cdots \le \g_{-2} \le \g_{-1} \le 0 \le \g_{1} \le
\g_2 \le \cdots$$ be the imaginary parts of the
zeroes of $L(s,\pi)$.  Let
$g$ and $h=\hat{g}$ be chosen so that $h\ge 0$
on $\R$ and let $c>0$ be a
cutoff parameter.  Then the function
$h_m(r)=h(r)(c^2-r^2)$
 is positive exactly when
$|r|<c$ and is the Fourier transform of $g_m=c^2g+g''$.
The support of $g_m$ is
of course also contained in $[-\log 2,\log 2]$
provided $g$ is suitably
regular.  If the sum
$2\Re\sum_{j=1}^m l_m(\eta_j) + g_m(0)\log D$ in
(\ref{simplexplic}) is
positive, then $\gamma_1 < c$ or $\gamma_{-1} > -c$,
i.e. $L(s,\pi)$
has a small zero.          To summarize:

\subsection{Criteria}\label{maincrit}

Our strategy will then be, for given archimedean parameters
$\eta_j$ and conductor $D$,
to find a function $g$ of support contained in $[-\log 2,\log 2]$
and for which
either
$$2\Re\sum_{j=1}^m l(\eta_j) + g(0)\log D < 0 $$
(which shows the L-function does not exist) or

$$2\Re\sum_{j=1}^m l_m(\eta_j) + g_m(0)\log D > 0$$
(which shows that it must have a {\em low zero} or
violate GRH in an effective range, as discussed below).

\subsection{What low zeroes mean without GRH}\label{wogrh}

Even if we do not assume GRH, we may still conclude from
$$2\Re\sum_{j=1}^m l_m(\eta_j) + g_m(0)\log D > 0$$
that the sum
$$\sum h_m(\g) > 0.$$  Thus, there are zeroes $\rho=\f12+i\g$ in the
region where $h_m(\g)>0$.
  We can explicitly compute the functions $h_m$ for our
choices of $g$ and examine where they are positive and negative within
the critical strip.  Since the density of zeroes increases
only logarithmically with their height (with an effective constant),
 and our functions $h_m(z)$ decay polynomially
 as $z\rightarrow
\infty$ in the critical strip, the zero must be
contained in an effectively bounded
region of the critical strip.

\begin{figure}
\hspace{-.65in}{\epsfysize=.6in
\epsfbox{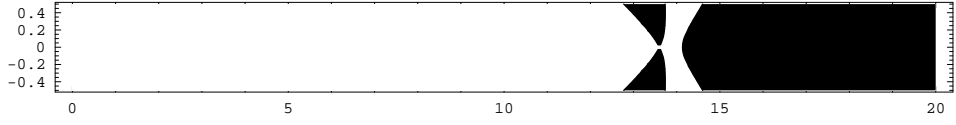}}
\hfill
\caption[Contour plot of $\Re h_{1m}(x+iy)$]{A
 contour plot of the function $\Re h_{1m}(x+iy)$.
We have colored the positive set white and the negative one
 black.\label{h1msupport}}
\end{figure}

As an example, \figref{h1msupport} is a contour plot
 of the function $h_{1m}$ defined
at the end of \secref{library}.  The white regions are where
$\Re h_{1m} > 0$, the black where $\Re h_{1m}<0$.

Figures~\ref{h2msupport} and \ref{h3msupport} contain plots for the
other functions we use.

\section{A library of functions}\label{library}

The main functions we use in this paper are
$$g_{1,p}(x)=\({\frac{\left( \pi
        \left( 1 - {\frac{|x|}{p}} \right)
        \cos ({\frac{\pi x}{p}}) \right)  +
     \sin ({\frac{\pi |x|}{p}})}{\pi }}\)/\cosh(x/2),$$

$$g_{2,p}(x)=\({\frac{4\pi \left( 1 - {\frac{|x|}{p}} \right)
        + 2\pi \left( 1 - {\frac{|x|}{p}}
         \right) \cos ({\frac{2\pi x}{p}}) +
      3\sin ({\frac{2\pi|x|}{p}})}{6\pi }}\)/\cosh(x/2),$$
and
$$g_{3,p}(x)=
\frac{\({\frac{54\pi \left( 1 - {\frac{|x|}{p}}
        \right) \cos ({\frac{\pi x}{p}}) +
     6\pi \left( 1 - {\frac{|x|}{p}} \right)
      \cos ({\frac{3\pi x}{p}}) +
     27\sin ({\frac{\pi |x|}{p}}) +
     11\sin ({\frac{3\pi |x|}{p}})}{60\pi }}\)}{\cosh(x/2)}.$$

We have normalized $g_{j,p}(0)=1$ and will often write
$g_j(x)=g_{j,\log(2)}(x)$. Ignoring the $\cosh(x/2)$'s
temporarily, the functions $g_{j,p}$ are rescalings of the
convolutions of $(\cos(\f{\pi x}{2}))^{j}$ with itself. Without
the $\cosh(x/2)$ term they would thus have a positive Fourier
transform on the real line, and the $\cosh(x/2)$ term spreads the
positivity into the critical strip. Were we to assume GRH we would
not need it.

\begin{lem} If an even function
$g(x)$'s Fourier transform is positive on the real line,
 then the Fourier transform of $g(x)/\cosh(x/2)$
is positive in the strip $-\f12 < \Im r < \f12$.
\end{lem}
{\bf Proof}: The Fourier transform of $\sech(x)$ is
$$\int_{-\infty}^\infty \f{2}{e^x+e^{-x}}e^{irx}dx~~~ =~~~
\pi\sech(\pi r/2).$$  This has positive
real part for $-\f12 < \Im r <\f12$ and the Fourier transform
converts multiplication into convolution, so the smeared
$\widehat{g\cdot\sech}$ remains positive in this strip.
\bx

We defined modified functions
 $$g_m=\f{c^2g(x)+g''(x)}{c^2g(0)+g''(0)},$$ which also have
$g_m(0)=1$.  (Of course we multiplicatively normalize $g(0)=g_m(0)=1$
to compare the explicit formulas from various test functions.)
These are used for showing the presence of low zeroes,
and since we do not assume GRH for this, we will actually use
$g=g_{1,p}(x)\cosh(x/2), g_{2,p}(x)\cosh(x/2),$ or
$g_{3,p}(x)\cosh(x/2)$.
Thus
$$g_{1m,p,c}(x) = {\frac{{\frac{{{\pi }^2}
         \left( -1 + {\frac{|x|}{p}} \right)
         \cos ({\frac{\pi x}{p}})}{{p^2}}} +
     {\frac{\pi \sin ({\frac{\pi |x|}{p}})}
       {{p^2}}} - {\frac{{c^2}
         \left( -\left( \pi
              \left( 1- {\frac{|x|}{p}} \right)
              \cos ({\frac{\pi x}{p}}) \right)
           - \sin ({\frac{\pi |x|}{p}}) \right) }
         {\pi }}}{{c^2} -
     {\frac{{{\pi }^2}}{{p^2}}}}}$$

$$g_{2m,p,c}(x) = {\frac{{\frac{{c^2}
         \left( 4\pi
            \left( 1- {\frac{|x|}{p}} \right)  +
           2\pi \left( 1- {\frac{x}{p}}
               \right)
            \cos ({\frac{2\pi x}{p}}) +
           3\sin ({\frac{2\pi |x|}{p}})
           \right) }{6\pi }} -
     {\frac{{\frac{8{{\pi }^3}
             \left( 1- {\frac{x}{p}} \right)
             \cos ({\frac{2\pi x}{p}})}{{p^2}}
           } + {\frac{4{{\pi }^2}
             \sin ({\frac{2\pi |x|}{p}})}{{p^2}}
           }}{6\pi }}}{{c^2} -
     {\frac{4{{\pi }^2}}{3{p^2}}}}}$$
and
$$g_{3m,p,c}(x)=~~~~~~~~~~~~~~~~~~~~~~~~~~~~~~~~~~
~~~~~~~~~~~~~~~~~~~~~~~~~~~~~~~~~~~~~~~$$
\begin{eqnarray}
  \frac{c^2\left( 54\pi
            \left( 1 - {\frac{|x|}{p}} \right)
            \cos ({\frac{\pi x}{p}}) +
           6\pi \left( 1 - {\frac{|x|}{p}} \right)
            \cos ({\frac{3\pi x}{p}}) +
           27\sin ({\frac{\pi |x|}{p}}) +
           11\sin ({\frac{3\pi |x|}{p}})
            \right) }{60\pi\(c^2 - \f{9\pi^2}{5p^2}\)}~~
+ & &\nonumber \\
  \frac{54\pi^3
             \left( -1 + {\frac{|x|}{p}} \right)
             \cos ({\frac{\pi x}{p}})  +
          54{{\pi }^3}
             \left( -1 + {\frac{|x|}{p}} \right)
             \cos ({\frac{3\pi x}{p}})
 + 81{{\pi }^2}
             \sin ({\frac{\pi |x|}{p}}) -
          63{{\pi }^2}
             \sin ({\frac{3\pi |x|}{p}})
            }{60\pi p^2\({c^2} -
     {\frac{9{{\pi }^2}}{5{p^2}}}\)}~. & & \nonumber
\end{eqnarray}
Since we are interested in finding zeroes in the range
from $\f 12 \pm 14.13472i$, we will now take $c=14.13472$
and write
$$g_{1m}(x)=g_{1m,\log(2),14.13472}(x),$$
$$g_{2m}(x)=g_{2m,\log(2),14.13472}(x),$$
and
$$g_{3m}(x)=g_{3m,\log(2),14.13472}(x).$$

The Fourier transforms of these functions are

$$h_{1m}(r)={\frac{-8{p^3}{{\pi }^2}
     \left( {c^2} - {x^2} \right)
     {{\cos ({\frac{px}{2}})}^2}}{\left( -{c^2}
          {p^2}   + {{\pi }^2} \right)
     {{\left( \pi  - px \right) }^2}
     {{\left( \pi  + px \right) }^2}}},$$

$$h_{2m}(r)={\frac{128p{{\pi }^4}\left( {c^2} - {x^2} \right)
     {{\sin ({\frac{px}{2}})}^2}}{\left( 3{c^2}
        {p^2} - 4{{\pi }^2} \right)
     {{\left( -4{{\pi }^2}x + {p^2}{x^3} \right) }^
       2}}},$$
and
$$h_{3m}(r)={\frac{2304{p^3}{{\pi }^6}
     \left( {c^2} - {x^2} \right)
     {{\cos ({\frac{px}{2}})}^2}}{\left( 5{c^2}
        {p^2} - 9{{\pi }^2} \right)
     {{\left( 9{{\pi }^4} -
          10{p^2}{{\pi }^2}{x^2} + {p^4}{x^4}
           \right) }^2}}}.$$

\begin{figure}
\hspace{-.65in}{\epsfysize=.6in
\epsfbox{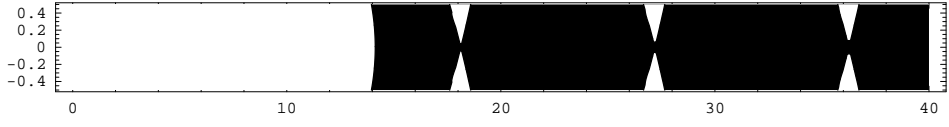}}
\hfill
\caption[Contour plot of $\Re h_{2m}(x+iy)$]{A
 contour plot of the function $\Re h_{2m}(x+iy)$.
We have colored the positive set white and the negative one
 black.\label{h2msupport}}
\vspace{1cm}
\end{figure}

\begin{figure}
\hspace{-.65in}{\epsfysize=.6in
\epsfbox{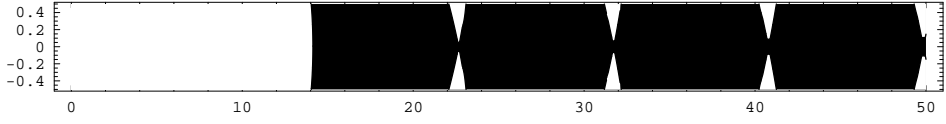}}
\hfill
\caption[Contour plot of $\Re h_{3m}(x+iy)$]{A
contour plot of the function $\Re h_{3m}(x+iy)$.
We have colored the positive set white and the negative set
 black.\label{h3msupport}}
\end{figure}

We show the contour plots of the functions $h_{2m}$ and $h_{3m}$
in Figures~\ref{h2msupport} and \ref{h3msupport},
the plot of $h_{1m}$ having been presented above in \figref{h1msupport}.

\newpage

\section{The highest lowest zero for $\pi_\infty$ real}

We restate
\begin{thm}[=1.1]\nonumber
 Let $\pi=\otimes_{p\le\infty}\pi_p$
be a cuspidal automorphic representation of $GL_n$ over $\Q$
 with a trivial central
character and whose archimedean type $\pi_\infty$ is real.  Then
$L(s,\pi)$ has a {\em low zero}.
\end{thm}

First we will note that for a fixed degree $m$, L-functions with
large $\eta_j$'s or large conductor $D$ must have low zeroes. This
is because Stirling's formula implies that
$$l(\eta)=\f{1}{2\pi} \int_{{\R}}
h(r)\f{\G_{{\R}}'}{\G_{{\R}}}\(\f{1}{2}+\eta+ir\)dr$$
has a positive real part
for $\eta$ large.  Thus, the lowest zero is only an issue
 for ``small'' archimedean parameters $\eta_j$ and
 small conductor --
partly because
$l(\eta)$ is bounded from below in
 $\Re\eta>-\f{1}{2}$ (which we may assume by (\ref{jsbound})).

We will present two different proofs of \thmref{realarch}.


%
%

\vspace{.5cm}

\begin{figure}
{\epsfysize=3in \epsfbox{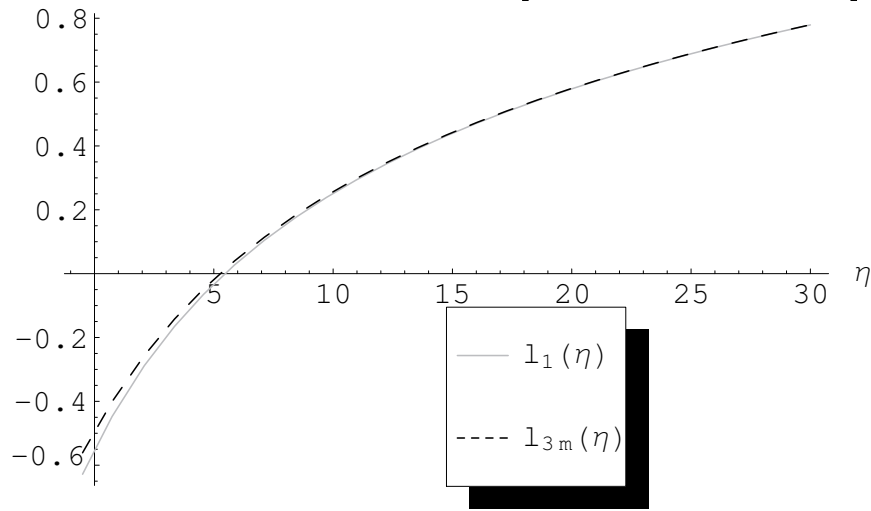}}
\vspace{1 cm}
\caption[The functions $l(\eta)$ and $l_{3m}(\eta)$]{
The functions $l(\eta)$ and
$l_{3m}(\eta)$.\label{functionsl1andl3m}}
\end{figure}

\begin{figure}
{\epsfysize=3in \epsfbox{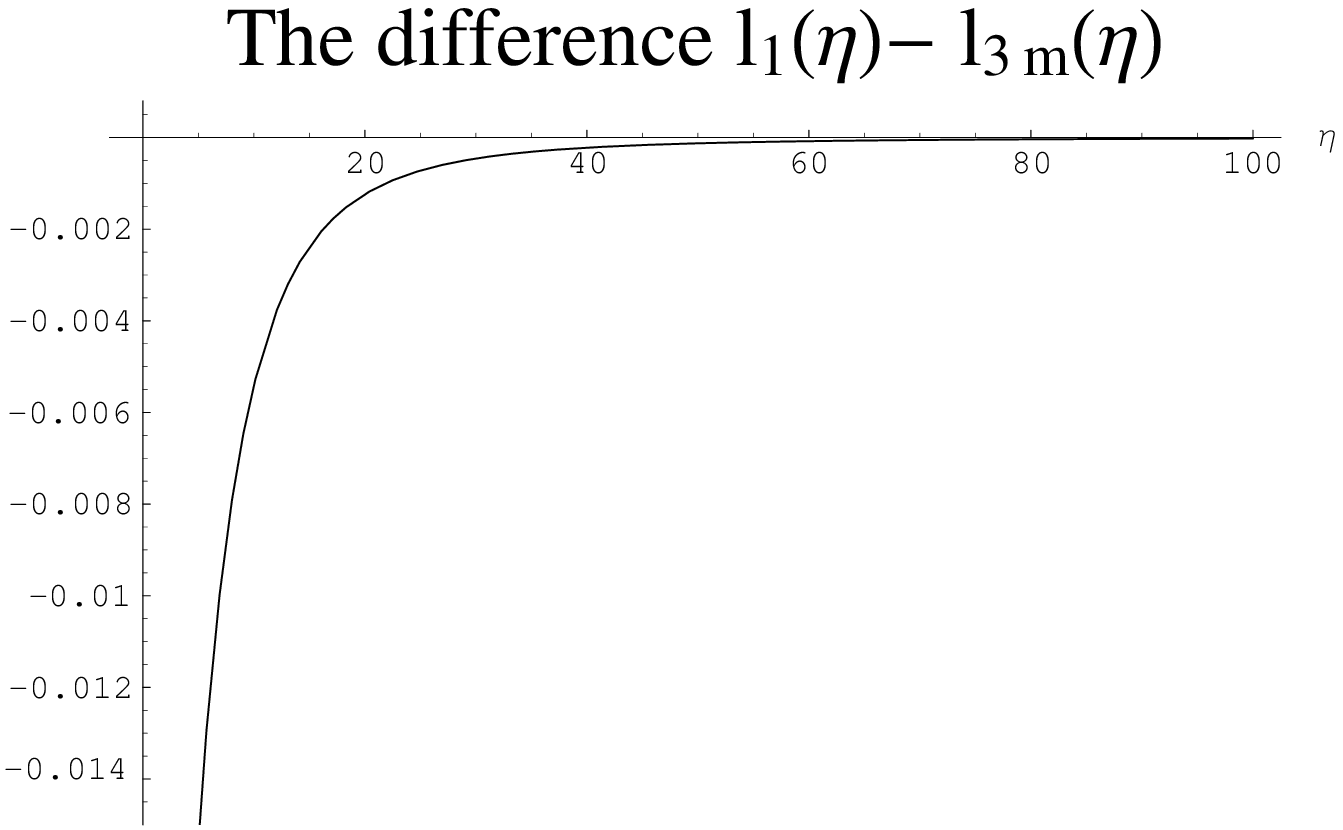}}
\vspace{.1 cm}
\caption{The difference between $l(\eta)$ and
$l_{3m}(\eta)$ \label{differencel1-l3m}}
\end{figure}

\begin{figure}
{\epsfysize=2in \epsfbox{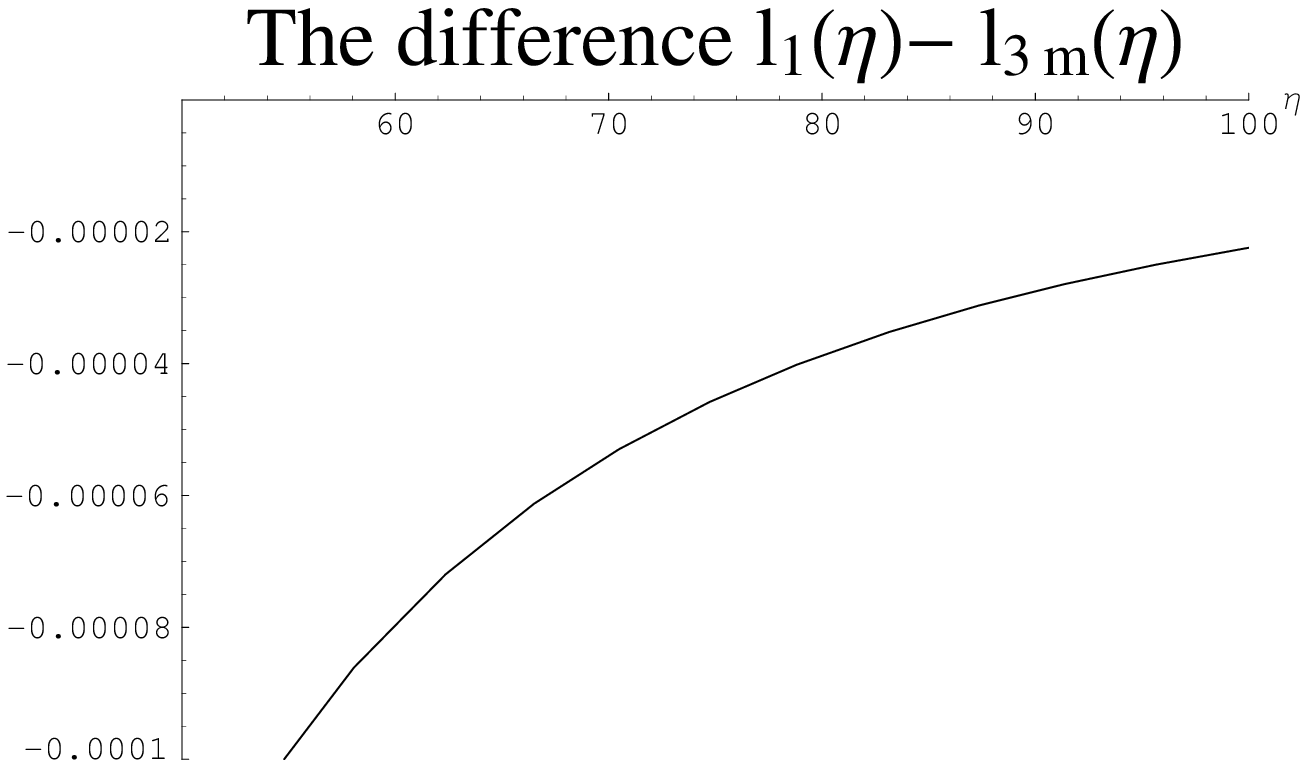}} \vspace{1 cm}
\caption{The difference between $l(\eta)$ and $l_{3m}(\eta)$,
magnified.\label{differencel1-l3mzoom}}
\end{figure}

{\bf Picture Proof of \thmref{realarch}}:
Figures~\ref{functionsl1andl3m},~\ref{differencel1-l3m}, and
 \ref{differencel1-l3mzoom}
indicate that $l_1(\eta)< l_{3m}(\eta)$ for $\eta\ge -\f12$, so
the theorem follows from Criteria~\ref{maincrit}.

 \bx

{\bf Less-Pictorial Proof of \thmref{realarch}}: This proof also
relies on numerical computation, but demonstrates how a proof can
be made even if the function $l$ is not strictly less than the
modified $l_{m}$. It uses $l_{2m}$ instead of $l_{3m}$.

We noted before in Criteria \ref{maincrit} that if
$$2\sum_{j=1}^m l_{2m}(\eta_j)+\log D \ge 0$$ then there is a indeed
a {\em low zero}, while if
$$2\sum_{j=1}^m l_{1}(\eta_j)+\log D \le 0,$$ the L-function
actually cannot exist to begin with.  Thus we are reduced to
dismissing the situation where
$$\sum_{j=1}^m l_{2m}(\eta_j) < 0$$ and
$$\sum_{j=1}^m\left[l_1(\eta_j)-l_{2m}(\eta_j)  \right] > 0$$
hold simultaneously.
Partition the $\eta_j \in (-\f12,\infty)$ into 3 sets:

$$N=\left\{ \eta_j \mid l_{2m}(\eta_j) \le 0,
l_1(\eta_j)-l_{2m}(\eta_j)\le 0
\right\}~~=~~(-\f12~,~5.4471\cdots],$$

$$S=\left\{ \eta_j \mid l_{2m}(\eta_j) >  0,
l_1(\eta_j)-l_{2m}(\eta_j)\le 0
\right\}~~=~~(5.4472\cdots~,~8.6553\cdots],$$
and
$$P=\left\{ \eta_j \mid l_{2m}(\eta_j) > 0,
l_1(\eta_j)-l_{2m}(\eta_j) > 0
\right\} ~~=~~(8.6553\cdots~,~\infty).$$
Of course if $$\sum_{j=1}^m l_{2m}(\eta_j)<0,$$ then also
$$\sum_{\eta_j\in N\cup P} l_{2m}(\eta_j)<0,$$ and if
$$\sum_{j=1}^m \left[l_1(\eta_j)-l_{2m}(\eta_j)\right] > 0,$$ then
$$\sum_{\eta_j\in N\cup P} \left[l_1(\eta_j)-l_{2m}(\eta_j)\right]
 > 0$$ as well.  Thus we need only consider the case where $S$ is
empty.

\begin{figure}
{\epsfysize=3in
\epsfbox{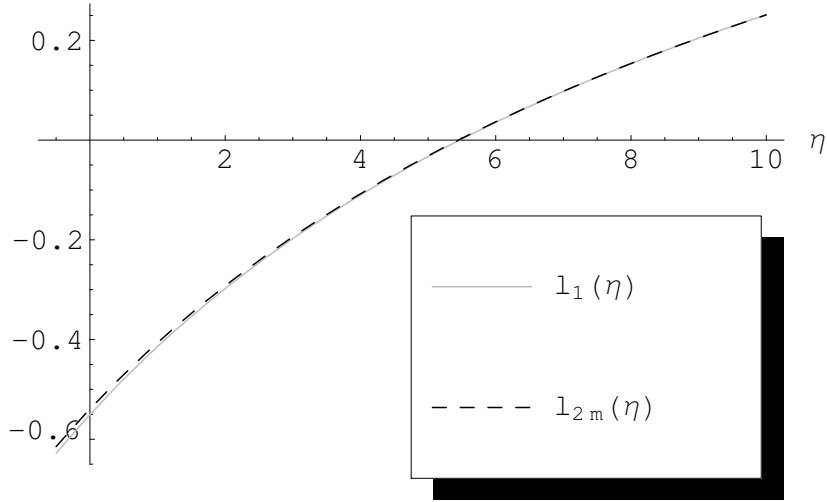}}
\caption{The graphs of $l(\eta)$ and
$l_{2m}(\eta)$.\label{functionsl1andl2m}}
\end{figure}

From computer investigations (see \figref{functionsl1andl2m}) on
the functions $l_1(\eta)$ and $l_{2m}(\eta)$ we can determine the
following very precise
 information:

\begin{figure}
{\epsfysize=3in \epsfbox{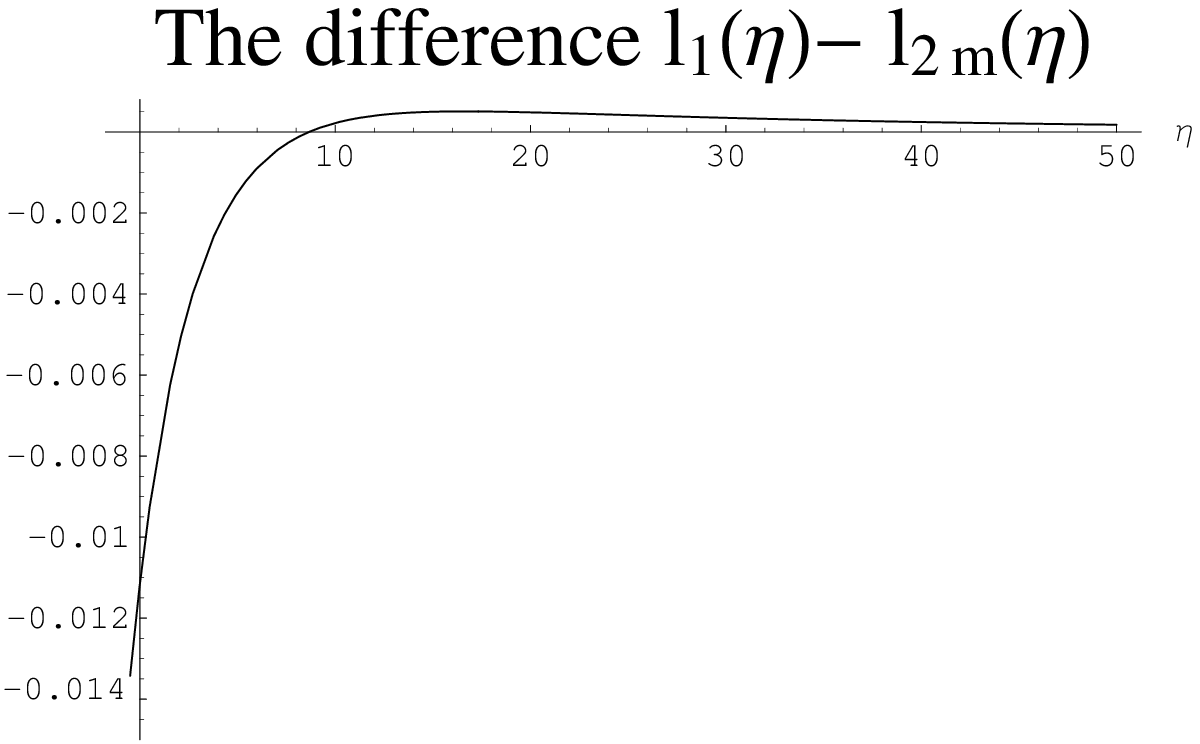}}
\caption[The
 difference between $l(\eta)$ and
$l_{2m}(\eta)$]{ The
 difference between $l(\eta)$ and
$l_{2m}(\eta)$.\label{differencel1-l2m}}
\end{figure}

$$\eta_j \in N \Longrightarrow -.628291 \le l_1(\eta_j) \le
0~~~~,~~~~l_1(\eta_j)-l_{2m}(\eta_j) \le -.001201$$
and
$$\eta_j \in P \Longrightarrow  l_1(\eta_j) \ge .187484~~~~~~~,~~~~~~~
0 \le l_1(\eta_j)-l_{2m}(\eta_j) \le .0005801.$$

Thus
$$0 > \sum_{j=1}^m l_{2m}(\eta_j) = \sum_{\eta_j\in N} l_{2m}(\eta_j)
+ \sum_{\eta_j\in P} l_{2m}(\eta_j)$$
$$\ge |N|(-.628291)+|P|(.187484),$$
which implies
$$\f{|N|}{|P|} < \f{.187484}{.628291}= .298403.  $$

On the other hand
$$0<\sum_{j=1}^m\left[l_1(\eta_j)-l_{2m}(\eta_j)    \right]
=\sum_{\eta_j\in N}  \left[l_1(\eta_j)-l_{2m}(\eta_j)   \right]
+\sum_{\eta_j\in P}  \left[l_1(\eta_j)-l_{2m}(\eta_j)   \right]
$$
$$\le |N|(-.001201)+|P|(.005801)$$
forces
$$\f{|N|}{|P|} > \f{.0005801}{.001201} = .483104,$$ a contradiction.
\bx

\section{Low zeroes for modular form L-functions}

In this section we prove that L-functions of cusp forms
on $GL_2$ over $\Q$ have low zeroes:

\begin{thm}[=1.2] Let $\pi$ be a cuspidal automorphic representation
of $GL_2$ over $\Q$ with a trivial central character.  Then
$L(s,\pi)$ has a {\em low zero}.
\label{gl2}
\end{thm}

Before giving the proof we shall give some background on the
hardest case -- Maass form L-functions.  In particular we will
precisely describe their completions, analytic continuations,
and functional equations in some important cases.

\subsection{Background on Maass forms on $\G_0(p)\backslash \U$}
It is known that if $D=p$ is a prime and $\pi$ is a cuspidal
automorphic representation not corresponding to a holomorphic
modular form, then $\pi$ instead corresponds to a Maass form
$\phi$ on $\G_0(p)\backslash \U$.  The Laplace operator $\D$,
the Hecke operators $T_n,n\ge 0$, as well as the involutions
$$(T_{-1}f)(x+iy)=f(-x+iy)$$
$$(W_pf)(x+iy)=f\(\f{-1}{p(x+iy)}\)$$
all commute.  Thus, after diagonalizing, we may take a basis
of Maass cusp forms on $\G_0(p)\backslash \U$ which are joint
eigenfunctions of $\D,T_n,T_{-1},$ and $W_p$.  Writing
$$\D\phi=\l\phi~~,~~\l=\f14-\nu^2,$$
$\phi$ has the Fourier expansion
$$\phi(x+iy)=\sum_{n\in{\Z}} c_n\sqrt{y}
K_{\nu}(2\pi|n|y)e^{2\pi inx},$$
where
$$K_{\nu}=K_{-\nu}=\f12 \int_0^\infty
e^{-y(t+t\i)/2}t^\nu \f{dt}{t}$$ is the $K$-Bessel function
of order $\nu$.  The cuspidality condition forces $c_0=0$;
the involution $T_{-1}$ interchanges $c_n$ and $c_{-n}$.

There are four symmetry classes of Maass forms under the action
of the involutions $T_{-1}$ and $W_p$.  The standard argument of
Hecke and Maass to prove that the L-functions of cusp forms are
entire also describes the functional equations
of L-functions of Maass forms having various symmetries.

\begin{prop}
Suppose $\phi$ is a Maass form on $\G_0(p)\backslash \U$
with
$$T_{-1}\phi=(-1)^\tau\phi$$ and
$$W_p\phi=(-1)^\omega\phi~~,~~\tau,\omega = 0\mbox{ or }1.$$
Multiplicatively normalize the coefficients of $\phi$ so that
$a_1=1$ and
$$\phi(x+iy) = \left\{  \begin{array}{ll}
\sum_{n=1}^\infty a_n\sqrt{y}K_{\nu}(2\pi ny)\cos(2\pi n x), &
\tau=0 \\
\sum_{n=1}^\infty a_n\sqrt{y}K_{\nu}(2\pi ny)\sin(2\pi n x), &
\tau = 1.\\
\end{array}
    \right.$$

Then $$\L(s,\phi)=\G_{{\R}}(s+\tau+\nu)
 \G_{{\R}}(s+\tau-\nu)\sum_{n=1}^\infty \f{a_n}{n^s}$$
satisfies the functional equation
\begin{equation}
\L(s,\phi)=(-1)^{\tau+\omega}p^{1/2-s}\L(1-s,\phi).
\label{maassfunceq}
\end{equation}
\label{maassl}
\end{prop}

{\bf Proof}: First consider the case $\tau=0$.  Then
$$\int_0^\infty \phi(iy)y^{s-1/2}\f{dy}{y} =\
\sum_{n=1}^\infty a_n \int_0^\infty K_{\nu}(2\pi ny)y^s\f{dy}{y}$$
$$=\sum_{n=1}^\infty a_n(2\pi n)^{-s}\left[
\int_0^\infty K_\nu(y)y^s\f{dy}{y}      \right]$$
$$=\sum_{n=1}^\infty a_n(2\pi n)^{-s}\left[2^{s-2}\G(\f{s+\nu}{2})
\G(\f{s-\nu}{2})\right]$$
$$=\f14 \L(s,\phi).$$
The transformation property
$$\phi(iy)=(-1)^{\omega}\phi(\f{i}{py})$$
gives
$$\L(s,\phi)=4\int_0^\infty \phi(iy)y^{s-1/2}\f{dy}{y}
= 4(-1)^\omega \int_0^\infty \phi(\f{i}{py})y^{s-1/2}\f{dy}{y}$$
$$=4(-1)^\omega \int_0^\infty \phi(\f{iy}{p})y^{1/2-s}\f{dy}{y}$$
$$=4p^{1/2-s}(-1)^\omega \int_0^\infty \phi(iy)y^{1/2-s}\f{dy}{y}$$
$$=(-1)^\omega p^{1/2-s}\L(1-s,\phi).$$

If instead $\tau=1$ then actually $\phi(iy)=0$ and
we instead consider the derivative
$$\phi'(x+iy):=\f{d}{dx}\phi(x+iy)$$
$$=\sum_{n=1}^\infty (2\pi n)a_n\sqrt{y}
K_{\nu}(2\pi ny)\cos(2\pi nx).$$
The action under $W_p$ now reads
$$\phi'(iy)=(-1)^\omega\phi'(\f{i}{py^2})(\f{-1}{py^2}).$$
We also have that
$$\int_0^\infty \phi'(iy)y^{s+1/2}\f{dy}{y}
= \sum_{n=1}^\infty (2\pi n)a_n \int_0^\infty K_\nu(2\pi ny)
y^{s+1}\f{dy}{y}$$
$$=\sum_{n=1}^\infty a_n(2\pi n)^{-s}\left[2^{s-1}\G(\f{s+1+\nu}{2})
\G(\f{s+1-\nu}{2})\right],$$
and the functional equation for $\L(s,\phi)$ follows as before.\bx

\subsection{Low zeroes for Maass form L-functions}

We will first prove \thmref{gl2} for Maass forms through
a series of propositions.

\begin{prop}
\label{lowhighlevel}
Every Maass form L-function whose conductor satisfies
$$D\ge 3~~,~~\mbox{if }~~ T_{-1}\phi=\phi$$
or
$$D \ge 2~~,~~\mbox{if  }~~ T_{-1}\phi=-\phi$$
has a {\em low zero}.
\end{prop}
{\bf Proof:} In these two symmetry classes the gamma factors
of $\L(s,\phi)$ are either
$$\G_{{\R}}(s+\nu)\G_{{\R}}(s-\nu)$$
or
$$\G_{{\R}}(s+1+\nu)\G_{{\R}}(s+1-\nu),$$
depending on whether $\phi$ is even or odd under $T_{-1}$.
In each case we may assume the parameter $\nu$ is not real and hence
purely imaginary, because \thmref{realarch} already covers
the case of real archimedean type.

\begin{figure}
{\epsfysize=3in
\epsfbox{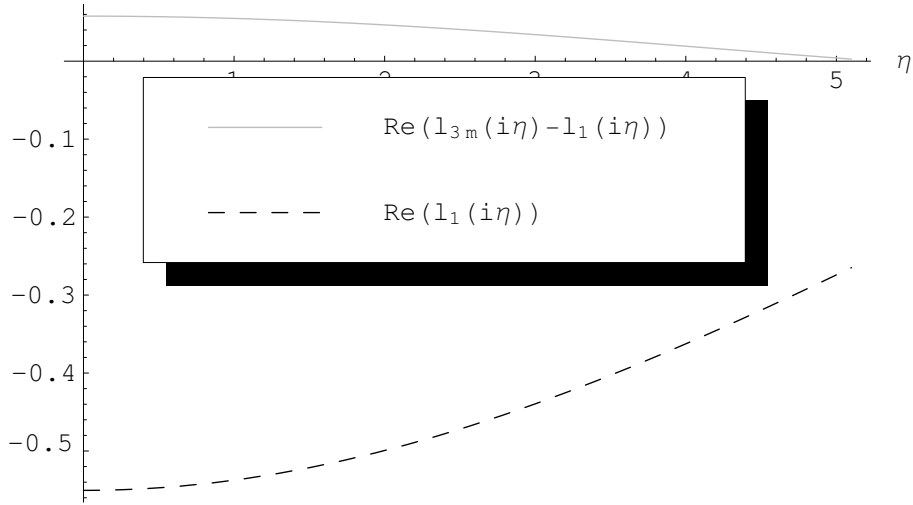}}
\caption[A plot of $l_{3m}(\eta)-l_1(\eta)$ and $l_1(\eta)$]{A plot
showing that
$\Re l_1(ir)+\f{\log 3}{4}<0$ and $\Re\(l_{3m}(ir)-l_1(ir)\) > 0$
 for $-5.1 < r < 5.1$.\label{l3m-l1andl1}}
\end{figure}

In the first case we have that
$$\Re \left( l_{3m}(ir)-l_1(ir) \right) >
 0~~~~\mbox{if } -5.1 < r <5.1,$$
a range in which $\Re l_1(ir)< -\f{\log 3}{4}\approx -0.274653$
(see \figref{l3m-l1andl1}).

In the second case
$$\Re\(l_{3m}(1+ir)-l_1(1+ir)\)>0~~~~\mbox{if } -5.5 < r <5.5,$$

\begin{figure}
{\epsfysize=3in \epsfbox{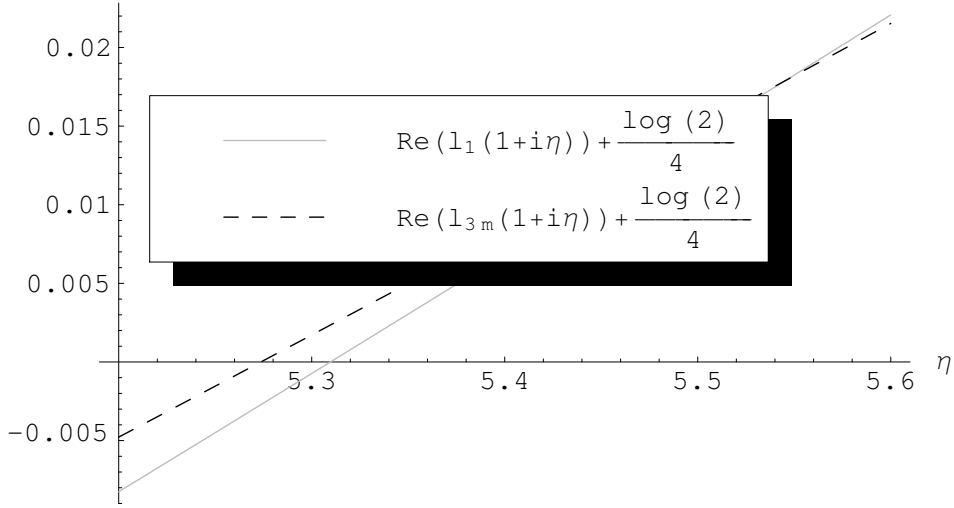}} \caption[A plot of
$l_{1}(\eta)$ and $l_{3m}(\eta)$]{A plot showing that $\Re
l_1(1+ir)+\f{\log 2}{4}<0$ and $\Re\(l_{3m}(1+ir)-l_1(1+ir)\) > 0$
 for $-5.5 < r < 5.5$.\label{l1andl3m}}
\end{figure}

where $\Re l_1(1+ir)$ and $\Re l_{3m}(1+ir)$ are both less than
$-\f{\log 2}{4}$ (see \figref{l1andl3m}).  Criteria \ref{maincrit}
thus shows there are {\em low zeroes} in either case.\bx

\vspace{.4 cm}

This next proposition handles the case of Maass forms at full
level (i.e. unramified for all primes $p<\infty$):

\begin{prop} If $\phi$ is a Maass form on
$SL_2({\Z})\backslash\U$ then $L(s,\phi)$ has a {\em low zero}.
\label{fulllevellow}
\end{prop}
{\bf Proof:} We again break the proof up into two cases, according
to whether $\phi$ is even or odd under $T_{-1}$.  By \thmref{realarch}
 we need only consider the case
$\Re\nu=0$.

\begin{figure}
{\epsfysize=3in
\epsfbox{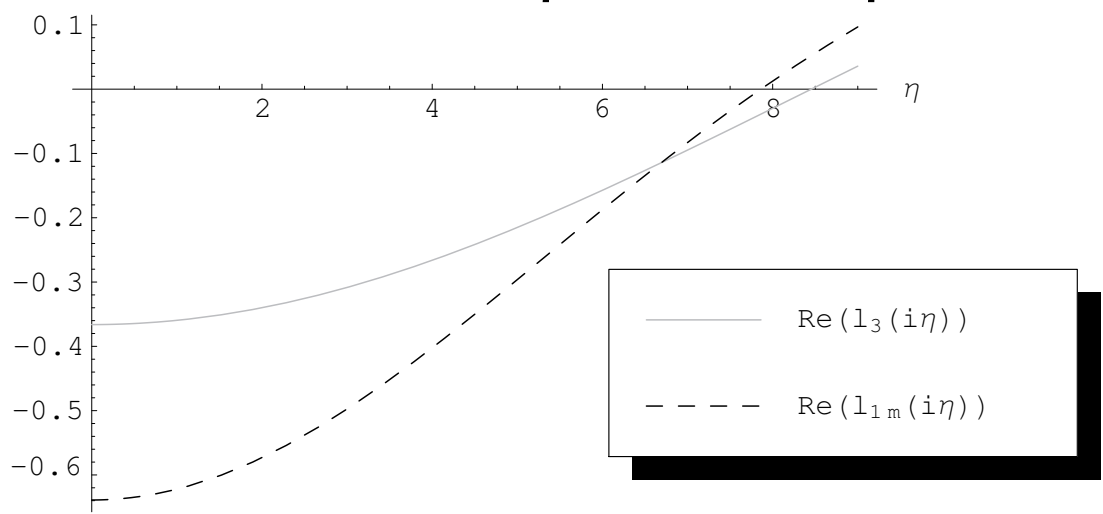}}
\caption{A plot showing that
$\Re l_3(\nu)$ is negative when $\Re l_{1m}(\nu)$ is.\label{l3andl1m}}
\end{figure}

If $\phi$ is even then the gamma factors of $L(s,\phi)$ are
$$\G_{{\R}}(s+\nu)\G_{{\R}}(s-\nu).$$  \figref{l3andl1m}
shows
$\Re l_3(\nu)$ is negative when $\Re l_{1m}(\nu)$ is,
which by Criteria \ref{maincrit} proves the proposition in this
case.

\begin{figure}
{\epsfysize=3in
\epsfbox{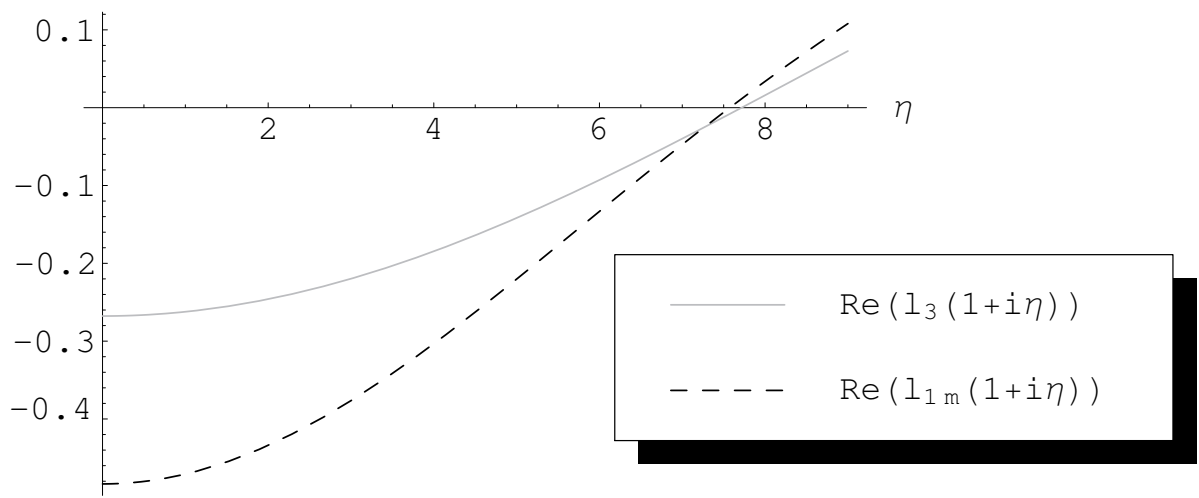}}
\caption{ A plot showing that
$\Re l_3(1+\nu)$ is negative when $\Re l_{1m}(1+\nu)$ is.
\label{l3andl1modd}}
\end{figure}

If instead $\phi$ is odd the gamma factors are instead
$$\G_{{\R}}(s+1+\nu)\G_{{\R}}(s+1-\nu),$$
and similarly $\Re l_{3}(1+\nu)$ is negative when $\Re
l_{1m}(1+\nu)$ is -- see \figref{l3andl1modd}.  The proposition
follows by invoking Criteria \ref{maincrit}.\bx

\vspace{.4cm}

To handle the remaining case, of even Maass forms on
$\G_0(2)\backslash\U$, we will use a result about the smallest
even eigenvalue of the laplacian there.  Perhaps
\propref{leveltwo} below can be proven without such explicit
information.

\begin{prop}
If $\phi$ is a Maass form on $\G_0(2)\backslash\U$
which is even under both $T_{-1}$ and $W_2$, then its
Laplace eigenvalue exceeds $\f14+6.14^2$.
\label{gtwoeigbd}
\end{prop}

Hejhal \cite{hejhal} has numerically computed that the first such
eigenvalue is $\approx \f14+8.922^2$.  We present the following
argument to demonstrate a technique.

\begin{figure}
{\epsfysize=2in
\epsfbox{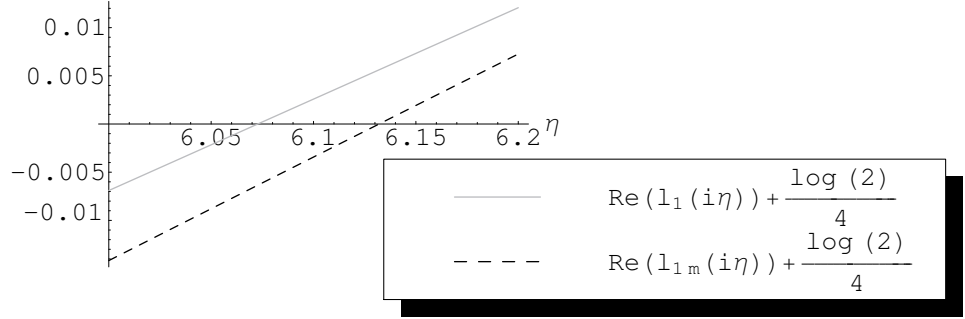}}
\caption[A plot of $l_{1m}(\eta)$ and $l_1(\eta)$]{ A plot showing that
 $\Re l_{1m}(ir)>-\f{\log 2}{4}$ for $r>6.135$ ,
while  $\Re l_1(ir)<-\f{\log 2}{4}$ for $r<6.07$\label{l1andl1m}}
\end{figure}

{\bf Proof:} First, \figref{l1andl1m} shows that

$$\Re{l_1(ir)+\f{\log 2}{4}}<0$$
for $-6.07 \le r \le 6.07$,
so we need only
consider the range $6.07 \le r \le 6.14$.  Using the symmetries,
set
$$f(t)=\phi\(\f{i}{\sqrt{2}}e^t\) = \sum_{n=1}^\infty
\f{a_n}{\sqrt{n}}W_{ir}\(\f{n}{\sqrt{2}}e^t\) = f(-t),$$
$$W_{ir}(y)=\sqrt{y} K_{ir}(2\pi y).$$
Thus $f$ is an even function in $t$, and so
\begin{equation}
f'(0)=\sum_{n=1}^\infty \f{a_n}{\sqrt{n}} W'_{ir}\(\f{n}{\sqrt{2}}\)
\f{n}{\sqrt{2}} = 0
\label{fpiszero}
\end{equation}
and
$$f'''(0)= \sum_{n=1}^\infty \f{a_n}{\sqrt{n}} \left[
W_{ir}'''\(\f{n}{\sqrt{2}}\)\f{n^3}{2\sqrt{2}} +
W_{ir}''\(\f{n}{\sqrt{2}}\)\f{3n^2}{2} +
W_{ir}'\(\f{n}{\sqrt{2}}\)\f{n}{\sqrt{2}}
\right]$$
\begin{equation}
=:\sum_{n=1}^\infty \f{a_n}{\sqrt{n}} V_{ir}(n) = 0.
\label{fpppiszero}
\end{equation}
The terms in the Fourier expansion decay rapidly with $n$,
and so we will use the first three terms as an approximation.
Recall that we are focusing on the range $6.07 \le r \le 6.14$.
We may assume that $\phi$ is a Hecke eigenform with $a_1=1$,
 and \cite{BDHI} have
proven that their coefficients satisfy the bound
$$|a_n| \le \tau(n)n^{5/28},$$
where $\tau(n)$ is the number of divisors of $n$.
Using the crude bound $\tau(n)\le 2 \sqrt{n}$ we can bound
the tails
\begin{equation}
\left|\sum_{n=4}^\infty
\f{a_n}{\sqrt{n}}W_{ir}'\(\f{n}{\sqrt{2}}\)\f{n}{\sqrt{2}}
        \right| \le 1.14\cdot10^{-7}
\label{fptail}
\end{equation}
and
\begin{equation}
\left|
\sum_{n=4}^\infty \f{a_n}{\sqrt{n}}V_{ir}(n)
\right| \le 2.7\cdot 10^{-5}.
\label{fppptail}
\end{equation}

Thus, (\ref{fpiszero}) and (\ref{fptail})
show
$$\left|
W_{ir}'\(\f{1}{\sqrt{2}}\)\f{1}{\sqrt{2}} + \f{a_2}{\sqrt{2}}
W_{ir}'\(\f{2}{\sqrt{2}}\)\f{2}{\sqrt{2}} +\f{a_3}{\sqrt{2}}
W_{ir}'\(\f{3}{\sqrt{2}}\)\f{3}{\sqrt{2}}
\right| \le 1.14\cdot 10^{-7}$$

while
(\ref{fpppiszero}) and (\ref{fppptail}) show
$$\left|
V_{ir}(\f{1}{\sqrt{2}}) + \f{a_2}{\sqrt{2}} V_{ir}(\f{2}{\sqrt{2}})
+ \f{a_3}{\sqrt{2}} V_{ir}(\f{3}{\sqrt{2}})
\right| \le 2.7\cdot 10^{-5}.$$

Now, $|W_{ir}'(\f{3}{\sqrt{2}})| \le 2.5\cdot 10^{-6}$
in the range $6.07 \le r \le 6.14$.  Yet
$W_{ir}(\f{1}{\sqrt{2}})$ and
 $W_{ir}(\f{2}{\sqrt{2}})\f{2}{\sqrt{2}}$
are much larger, never smaller than $5.7\cdot 10^{-5}$
in magnitude.  The ratio of
$$\f{ W_{ir}'( \f{1}{\sqrt{2}} )\f{1}{\sqrt{2}}}{W_{ir}'(
 \f{2}{\sqrt{2}} )\f{2}{\sqrt{2}}}$$
is smallest at $r=6.07$, where it is $\approx 1.475 > 1$.
Thus, we must have that $\f{a_2}{\sqrt{2}} > 1$ for
(\ref{fpiszero}) to be valid.

At the same time, such a value of $\f{a_2}{\sqrt{2}}$
is too large to achieve equality in (\ref{fpppiszero}).  This
is because it makes the second term much larger than
the first and third terms could possibly be with
the constraint that $a_3\le 2\cdot 3^{5/28}$:
$$|V_{ir}(\f{1}{\sqrt{2}})|+\f{2\cdot 3^{5/28}}{\sqrt{3}}
|V_{ir}(\f{3}{\sqrt{2}})| < |V_{ir}(\f{2}{\sqrt{2}})|~~,~~
6.07 \le r \le 6.14.$$

So (\ref{fpiszero}) and (\ref{fpppiszero}) cannot hold
simultaneously.  This contradiction shows every Maass form
on $\G_0(2)\backslash\U$ which is even under both $T_{-1}$
and $W_2$ has Laplace eigenvalue greater than $\f14+6.14^2$.
\bx

\begin{prop}
Maass form L-functions with conductor $D=2$
(which correspond to Maass forms on $\G_0(2)\backslash\U$)
 have {\em low zeroes}.
\label{leveltwo}
\end{prop}
{\bf Proof:} By \propref{lowhighlevel} we need only consider
the even Maass forms, where the gamma factors are
$$\G_{{\R}}(s+\nu)\G_{{\R}}(s-\nu).$$
In fact, by \propref{maassl} we can assume
that $\phi$ is even under both $W_2$ and $T_{-1}$;
otherwise (\ref{maassfunceq}) dictates
$$\L(\f12,\phi)=-\L(\f12,\phi)=0.$$
The function $\Re l_{1m}(ir)>-\f{\log 2}{4}$ for $r>6.135$
(\figref{l1andl1m}),
and \propref{gtwoeigbd} shows all even eigenvalues are in that
range.\bx

\vspace{.5cm}

{\bf Proof of \thmref{gl2}:} Every cuspidal automorphic
representation on $GL_2$ over $\Q$ comes from either a Maass form
or a holomorphic modular form. Both holomorphic modular forms and
non-tempered Maass forms (i.e. $\l<\f14$) have real archimedean
type and are thus covered under \thmref{realarch}.  The rest of
the Maass forms (the tempered ones) are covered by Propositions
\ref{lowhighlevel}, \ref{fulllevellow}, and \ref{leveltwo}. \bx

\section{Cuspidal eigenvalue bounds}

Now we move our focus completely towards automorphic
representations rather than on their L-functions. In this section
and in the next we will examine the discrete spectrum of the
laplacian $\D$ on $L^2(\quo n)$.\footnote{Of course our methods
carry over to some congruence covers but we will restrict our
attention to full-level here.} We normalize our laplacian so that
its continuous spectrum on $L^2(SL_n({\R})/SO_n({\R}))$ spans the
interval from
$$\l_1(SL_n({\R})/SO_n({\R})) = \f{n^3-n}{24}$$ to $\infty$.

Because the ring of invariant differential operators ${\cal R}$ on
$SL_n({\R})/SO_n({\R})$ is commutative, we may take a basis
of Laplace eigenfunctions which are also common eigenfunctions of the
operators in
${\cal R}$.  Thus, to each discrete eigenfunction $\phi\in
L^2(\quo n)$ we can attach Langlands parameters $\mu_1,\ldots,\mu_n$.
These describe $\phi$'s eigenvalues under the different operators
in $\cal R$; in particular, the Laplace eigenvalue satisfies
$$\D\phi=\l\phi~~,
~~\l=\f{n^3-n}{24}-\f{\mu_1^2+\cdots + \mu_n^2}{2}.$$
By the Jacquet-Shalika ``trivial'' bound \cite{JS}
$$|\Re{\mu_j}|< \f 12~~,~~j=1,\ldots,n.$$
Thus,
\begin{equation}
\l>\f 12 \sum_{j=1}^n(\Im\mu_j)^2+ \f{n^3-{\bf 4}n}{24}.
\label{jsspecbd}
\end{equation}
We will use (\ref{jsspecbd}) to bound $\l$ from below.

\subsection*{Positivity Functions}
Recall the function
$$g_{1,p}=\((1-\f{|x|}{p})\cos\(\f{\pi x}{p}\)+
\f{1}{\pi}\sin\(\f{\pi|x|}{p}\)\)/\cosh(x/2)~~,~~0<p\le \log 2.$$
Define
$$s(r)=\max \Re l_{1,\f{1}{2}}(ir+\sigma),$$
where the maximum is taken over $\sigma\in [-\f{1}{2},\f{3}{2}]$.

\begin{figure}
{\epsfysize=3in
\epsfbox{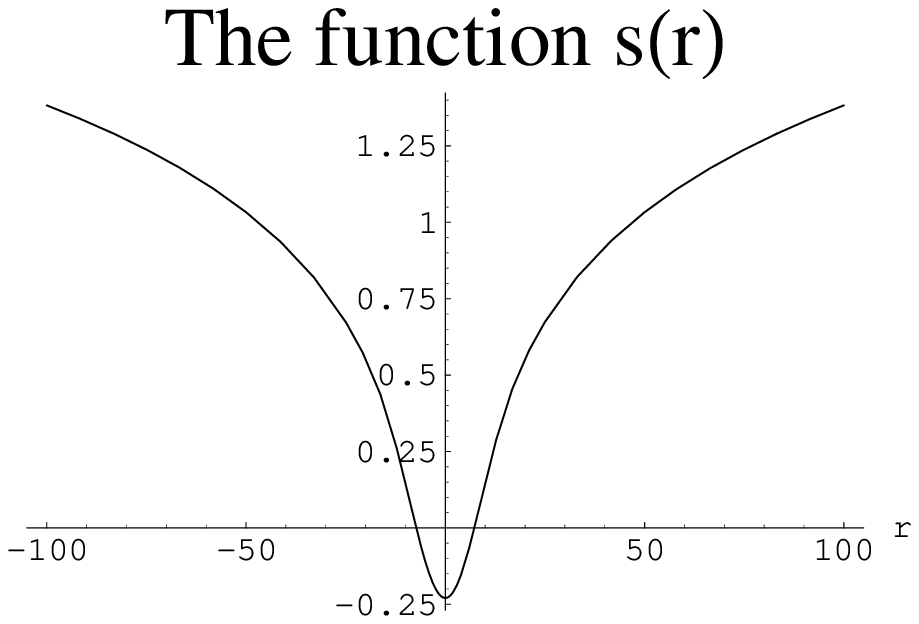}}
\caption{The graph of $s(r)$.\label{plotofs}}
\end{figure}

\subsection{Criteria}\label{secrit}
 If $\sum_{j=1}^n s(r_j)<0$ then there is no cuspidal
eigenfunction in $L^2(\quo n)$ whose Langlands parameters
$\mu_1,\ldots,\mu_n$ have $\Im\mu_j=r_j$.
\vspace{1cm}

If $\phi$ is a cusp form, then the {\em archimedean Ramanujan-Selberg
conjectures} assert that $\pi_\infty$ is tempered,
i.e. $\Re{\mu_j}=0$.  A consequence is that $\l^{cusp}
\ge \f{n^3-n}{24}$.  This was proven in \cite{Miller} unconditionally
using a similar positivity argument.  Here we can derive some stronger
results and different applications.

\begin{prop} (A trivial bound)
If $\phi$ is a cusp form in $L^2(\quo n)$, then with the
above notation
$$\sum_{j=1}^n r_j^2 > 51.84\(1+\f{1}{n-1}\).$$
\label{trivpos}
\end{prop}
{\bf Proof:} The plot shows $s(r)<0$ for $|r|<7.2$.  Thus
$\sum_{j=1}^n s(r_j)\ge 0$ only if at least one $|r_j|\ge 7.2.$
Since the $r_j$ are constrained to have $r_1+\cdots+r_n=0$, this means
$$\sum_ {j=1}^n r_j^2 \ge 7.2^2+(n-1)\(\f{7.2}{n-1}\)^2.$$\bx

\thmref{newlamb1} follows immediately from \propref{trivpos}
and (\ref{jsspecbd}).

\subsection{Extreme values} Given $d$ and the constraints
$\sum r_j^2=d,\sum r_j =0$, if the largest value obtained by $\sum
s(r_j)$ is negative, then Criteria \ref{secrit}
 implies $\l>\f{d}{2}+
\f{n^3-4n}{24}$.

\begin{prin}\footnote{Some may not consider the justification
to be a proof, but as we indicate, it can be verified
in the applications we use it for.}
If $(r_1,\ldots,r_n)$ is an extremal point of
$$\sum_{j=1}^n s(r_j)$$ subject to the constraints
$$\sum_{j=1}^n r_j=0~~,~~\sum_{j=1}^n r_j^2=d,$$
then the $r_j$ assume at most three distinct values.
\label{threepoints}
\end{prin}
{\bf Proof}: By Lagrange multipliers, there are real constants
$c_1,c_2\in\R$ such that
$$(s'(r_1),\ldots,s'(r_n))=c_1(r_1,\ldots,r_n)+c_2(1,\ldots,1),$$
i.e. the points $(r_j,s'(r_j))$ all lie on the intersection of some
line and the graph of $y=s'(x)$.

\begin{figure}
{\epsfysize=2.8in \epsfbox{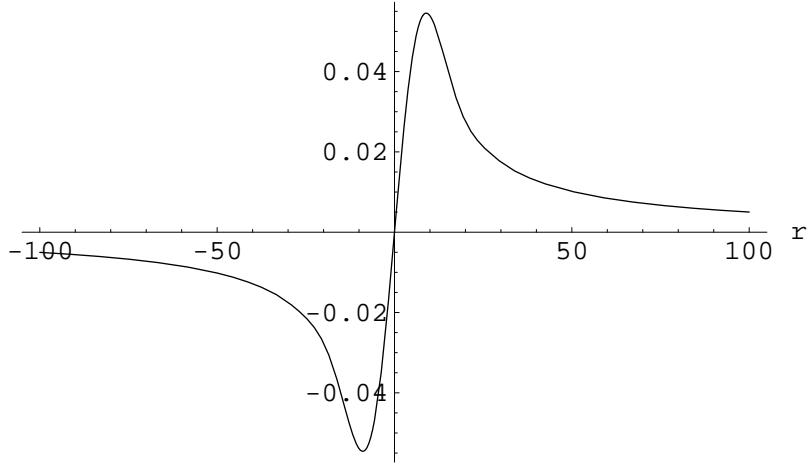}} \caption{ The graph
of $s'(r)$.\label{plotofsprime}}
\end{figure}

But no line crosses this graph in more than three places.
Even though \figref{plotofsprime} only shows the range $|x|\le 100$
it is legal to use this principle in this paper.  For another crossing
would give a value of $r_j$ so large that it would not enter into our
subsequent bounds.\bx

\begin{thm}We have the following bounds on the Laplace eigenvalue
of a cuspidal eigenfunction of $\D$ in $L^2(\quo n)$:

\begin{tabular}{||l||r|r|r|r|r|r||}
\hline\hline
$n$ & $3$ & $4$ & $5$ & $6$ & $7$ & $8$ \\
\hline
 $\l\ge$ & $87.625$ & $108.$& $140.875$ & $167$. & $201.125$ &
 $232.$
\\
\hline
    \end{tabular}

\label{lagbds}
\end{thm}

{\bf Proof}: By the last proposition, we need only
consider the case where there are $A r_1$'s, $B r_2$'s, and
$C r_3$'s, with
$r_1,r_2,r_3\in \R$,
$$Ar_1+Br_2+Cr_3=0~~,~~Ar_1^2+Br_2^2+Cr_3^2=d,$$ and then
try to find a large value of $d$ such that
$$As(r_1)+Bs(r_2)+Cs(r_3)$$ is always negative.  Since
$A,B,$ and $C$ are all positive integers which sum to $n$,
this is a finite calculation.  We will take $A,B,C>0$ by
allowing some of the values of $r_1,r_2$, and $r_3$ to coincide.
Then in terms of the parameter $r_3$, either
$$r_1= -{\frac{ACr_3 +
      {\sqrt{AB\left( - C^2
r_3^2   +
            A( D - C)
                + B( D  r_3^2- C r_3^2
                )  \right) }}}{A
      \left( A + B \right) }} $$
and
$$r_2= {\frac{-\left( BC r_3 \right)  +
     {\sqrt{AB\left( -C^2
              r_3^2 +
           A( D - C r_3^2)  +
           B( D - C r_3^2 )
            \right) }}}{B\left( A + B \right) }} $$
or instead
$$r_1= {\frac{-ACr_3 +
      {\sqrt{AB\left( -C^2 r_3^2 +
            A( D - C)
                + B( D  r_3^2- C r_3^2
                )  \right) }}}{A
      \left( A + B \right) }} $$
and
$$r_2= -{\frac{\left( BC r_3 \right)  +
     {\sqrt{AB\left( -C^2
              r_3^2 +
           A( D - C r_3^2)  +
           B( D - C r_3^2 )
            \right) }}}{B\left( A + B \right) }}.$$
Actually, the second set of solutions and the first are interchanged
upon $r_3\leftrightarrow -r_3$, so they take the same values.
For a given $n$, we need only enumerate the integer triples of
$A,B,C$ with $A\ge B\ge C>0$ and plot
$$As\(-{\frac{ACr_3 +
      {\sqrt{AB\left( - C^2
r_3^2   +
            A( D - C)
                + B( D  r_3^2- C r_3^2
                )  \right) }}}{A
      \left( A + B \right) }}
\) + $$
$$Bs\({\frac{-\left( BC r_3 \right)  +
     {\sqrt{AB\left( -C^2
              r_3^2 +
           A( D - C r_3^2)  +
           B( D - C r_3^2 )
            \right) }}}{B\left( A + B \right) }}
\)+Cs(r3)$$
 over the range
$$-{\frac{{\sqrt{A+ B}}{\sqrt{D}}}
    {{\sqrt{C}}{\sqrt{A + B + C}}}} \le r_3 \le
{\frac{{\sqrt{A + B}}{\sqrt{D}}}
    {{\sqrt{C}}{\sqrt{A + B + C}}}}.$$

 One finds the following values of $d$ work:

\begin{tabular}{||l||r|r|r|r|r|r||}
\hline\hline
$n$ & 3 & 4 & 5 & 6 & 7 & 8 \\
\hline
$d$ & 174 & 212 & 273 & 318 & 376 & 424 \\
\hline
    \end{tabular}

\bx

\begin{remark}
 \thmref{imrn} shows there exists a positive constant
$c>0$ such that
$$\l^{cusp} -\f{n^3-n}{24} > cn~~,~~n=1,2,\ldots.$$
The argument above
gives a much better constant.
\end{remark}

\subsection{Some open problems about $\l^{cusp}_1(\quo n)$}

\begin{conj} Fix $k=1,2,\ldots$ and denote the $k$-th cuspidal
eigenvalue of $\D$ on $L^2(\quo n)$ as $\l_k^{cusp}(\quo n)$.
 Then
the sequence
\begin{equation}
\left\{\f{\l_k^{cusp}(\quo n)-\f{n^3-n}{24}}{n}
\mid n=1,2,\ldots           \right\}
\label{eigdif}
\end{equation}
has a limiting distribution.
\end{conj}

\begin{ques}
Is the sequence in (\ref{eigdif}) also bounded from
above as well as from below?
\end{ques}

\section{Bounds on non-cuspidal eigenvalues}

Lubotzky asked if the bound $$\l\ge \f{n^3-n}{24}$$ could
also hold for the entire non-zero
discrete spectrum of $\D$ on
$L^2(\quo n)$, i.e. not just for cusp forms alone.
Although from the point of view of automorphic forms the cusp forms
are most essential, the entire discrete spectrum enters
into considerations in differential geometry.
In fact,
there are non-constant,
non-cuspidal, square-integrable residues of Eisenstein
series on $L^2(\quo n)$ which are discrete Laplace eigenfunctions,
and they are never tempered (that is, they violate $\Re\mu_j=0$).
The first example of one on $\quo n$ violating $\l \ge \f{n^3-n}{24}$
occurs for $n=68$:

\begin{thm}
There exists a discrete Laplace eigenfunction $$\phi\in
L^2(\quo {68})$$ such that
$$\D\phi=\l_\phi \phi~~,~~\l_\phi\approx 12916.6 < \f{68^3-68}{24}
=13098.5 .$$
Yet for $n\le 67$ the bound
$$\l \ge \f{n^3-n}{24}$$ is valid for every non-zero
discrete eigenvalue
of $\D$ on $L^2(\quo n)$.
\label{lubotz}
\end{thm}
Of course the failure of $\l \ge \f{n^3-n}{24}$ at
 $n=68$ is the typical case
for large $n$.

The key idea here is the classification of the discrete spectrum
in terms of cusp forms.  It was first conjectured by Jacquet
\cite{Jacquet} and later proven by M\oe glin-Waldspurger \cite{MoeWal}.
Let us now describe how discrete eigenfunctions can be constructed.
Factor $n=ra$ and let $\phi$ be a cusp form on $\quo a$.  The
group $SL_n({\R})$ has a rank $r-1$ parabolic subgroup $P$
of type $(a,a,\ldots,a)$ whose Levi component is
$$L=GL_a({\R})^r\cap SL_n({\R}).$$
The cusp form $\phi$ extends as a product
to the $r$ copies of $SL_a({\R})$
in $L$ in the obvious way.  Given $h=(h_1,\ldots,h_r)\in {\C}^r$
such that $h_1+\cdots+h_r=0$, we can form a character of the split
Levi component $A$ of $P$, and the Eisenstein series $E(P,g,\phi,h)$.
If $\phi$ has Langlands parameters $\mu_1,\ldots,\mu_a$, then
$E(P,g,\phi,h)$ has Langlands parameters
$$(\underbrace{\overbrace{\mu_1+h_1,\mu_2+h_1,\ldots,\mu_a+h_1}^a,\mu_1+h_2,\ldots\ldots\ldots,
\overbrace{\mu_1+h_r,\ldots,\mu_a+h_r}^a}_{ra=n}).$$  Furthermore,
$E(P,g,\phi,h)$ has a pole of order $r-1$ at
$h=(\f{r-1}{2},\f{r-3}{2},\ldots, -\f{r-1}{2})$ and its $r-1$st
iterated residue there is a discrete, $L^2$ eigenfunction of $\D$.
Moreover, all of them arise this way. For example, if $r=1$ these
are just cusp forms, and if $a=1$, constant functions.

We compute that the residue's Laplace eigenvalue is
$$2\l-\f{n^3-n}{12} =
-\sum_{j=1}^a \sum_{k=1}^r \(\mu_j+\f{r-1}{2}-k\)^2$$
$$=-r\sum_{j=1}^a \mu_j^2-a\f{r^3-r}{12}.$$

Incidentally, Maass forms with Laplace eigenvalue $\f14$ are known
to exist on congruence quotients of $SL_2({\R})/SO_2({\R})$. Using
this procedure one may already construct a discrete residue on a
congruence quotient of $SL_4({\R})/SO_4({\R})$ which violates the
$\l\ge \f{4^3-4}{24}$ bound.

\vspace{.5cm}

{\bf Proof of \thmref{lubotz}}:
Firstly, Hejhal (see \cite{hejhal}) has computed that
$\l_1^{cusp}(\quo 2) = 91.1413\cdots$, corresponding
to $\mu_1=-\mu_2\approx 9.534i$.  Thus the
Laplace eigenvalue of a residue formed from Hejhal's Maass form
has
$$2\l-\f{n^3-n}{12} \approx r(181.8)-2\f{r^3-r}{12},$$
and this difference is positive for
$$r < \sqrt{6\cdot 181.8+1} \approx 33.04.$$
For $r=34$ we have
$$\l\approx 12916.6 < \f{68^3-68}{24}=13098.5 .$$

If in fact there was an example of a residue for $n=ra<68$, with
$$-r\sum_{j=1}^a\mu_j^2-a\f{r^3-r}{12}<0,$$ we would necessarily
have
\begin{equation}
r > \sqrt{-\f{12}{a}\sum_{j=1}^a\mu_j^2+1}.
\label{rbdinlubo}
\end{equation}

We already know that $r>1$ since cusp forms obey the
$\l\ge\f{n^3-n}{24}$
bound.  Thus we can restrict to the cases $r\ge 2,a=1,\ldots,34$.
Using our pre-existing bounds for $\sum \mu_j^2$
we conclude $ra\ge 68$ -- see Table~\ref{lubopf} for details.

\begin{table}
\begin{tabular}{||r|r|r|r||}
\hline
\hline
$a$ & $r\le [\f{68}{a}]$ &
Lower bound for $-\sum_{j=1}^a\mu_j^2$ &
 Upper bound for $-\sum_{j=1}^a\mu_j^2$ \\
\hline
3           & 22 &  171.25& 120.75  \\
4 & 17 &    211          & 96.00  \\
5 & 13 &    271.75           & 70.  \\
6 & 11 &    316.5           & 60.  \\
7 &  9 &    374.25           &46.66  \\
8 & 8 &     422         & 42.00  \\
9 & 7 & 56.07 & 36.  \\
10 & 6 & 55.10 & 29.16  \\
11 & 6 & 54.27 & 32.08  \\
12 & 5 & 53.55 & 24.   \\
13 & 5 & 52.91 & 26.  \\
14 & 4  &  52.32 & 17.5  \\
15 & 4 & 51.79 &  18.75  \\
16 & 4 & 51.29 & 20.  \\
17 & 4 & 50.83 & 21.25  \\
18 & 3 & 50.38 &  12.  \\
19 & 3 & 49.96 & 12.66  \\
20 & 3 & 49.56 & 13.33  \\
21 & 3 & 49.18 & 14.  \\
22 & 3 & 48.80 & 14.66  \\
23 & 2 & 48.44 & 5.75  \\
24 & 2 & 48.09 & 6.  \\
25 & 2 & 47.75 & 6.25  \\
26 & 2 & 47.41 & 6.5  \\
27 & 2 &  47.08 & 6.75  \\
28 & 2 & 46.76 & 7.  \\
29 & 2 & 46.44 & 7.25  \\
30 & 2 &  46.12 & 7.5  \\
31 & 2 & 45.81 & 7.75  \\
32 & 2 & 45.51 & 8.  \\
33 & 2 & 45.20 & 8.25  \\
34 & 2 & 44.91 & 8.5  \\
\hline
\end{tabular}
\caption[Proof of \thmref{lubotz}]{This table completes the proof
of \thmref{lubotz}. Suppose a residue of a cusp form on $GL_a$
occurred on some $GL_n, n=ra<68$ with Laplace eigenvalue
$\le\f{n^3-n}{24}$.
  The second column gives the upper bound
$r \le [\f{68}{a}]$.  The third column gives a lower bound for
$-\sum_{j=1}^a\mu_j^2$ (from
 \propref{trivpos}),
 but the fourth
gives an upper bound on $-\sum_{j=1}^a\mu_j^2$ that would be
 satisfied by such a residue with low eigenvalue (as derived
in (\ref{rbdinlubo}) in the proof of \thmref{lubotz}).  The inconsistency
of these two
inequalities
is a contradiction which shows that the discrete Laplace spectrum
on $\quo n$ is contained in $\{0\}\cup[\f{n^3-n}{24},\infty)$
for $n<68$.}
\label{lubopf}
\end{table}

\bx

\section{Cuspidal cohomology}
 The positivity inequality can be applied to products of L-functions which
have poles, for example Rankin-Selberg L-functions
$L(s,\pi\otimes\tilde{\pi})$ of cuspidal automorphic forms $\pi$ on
$GL_n$.  If $\{\mu_{jk}\}_{j=1,k=1}^m$ are the archimedean $\G_{\R}$
parameters, the inequality reads
\begin{equation}
\int_{{\R}} g(x)\(e^{x/2}+e^{-x/2}\)dx +
2 \Re \sum_{j=1}^m \sum_{k=1}^m
l(\mu_{jk}) + g(0)\log D \ge 0.
\label{polepos}
\end{equation}
The new term in (\ref{polepos}) as compared to
(\ref{simplexplic}) comes from the
poles of $L(s,\pi\otimes\tilde{\pi})$. Also, here we
have  simply dropped the coefficients
entirely because $\f{L'}{L}(s,\phi)$ has
a Dirichlet series with
non-positive coefficients (see \cite{Rud-Sar} for
 a verification of this)
and so there is no restriction on the support of $g$.

If $\pi=\otimes_{p\le\infty}\pi_p$ comes from a
constant-coefficients cohomological cusp form on $GL_n(\A_{\Q})$
then $\pi_\infty$ is of either the form

$$\pi_{\infty}=\mbox{Ind}_{P_{(2,2,\ldots,2)}}^{GL_n}
(D_2,D_4,\ldots,D_n),~ n\mbox{ even},$$ or

$$\pi_{\infty}=\mbox{Ind}_{P_{(1,2,2,\ldots,2)}}^{GL_n}
(\hbox{sgn}(\cdot)^{\e},D_3,D_5,\ldots,D_n),~ n\mbox{ odd}.$$
($\hbox{sgn}$ is the sign character, $\e=0$ or 1, and $D_k$
denotes the $k$-th discrete series on $GL_2$, corresponding to
weight $k$ holomorphic forms.) Thus, if $n$ is written as $2m+t$,
$t=0$ or 1, the archimedean $\mu_{jk}$ can be computed via the
recipe summarized in \cite{Rud-Sar} and are the following
multisets:
\begin{eqnarray}
\{\mu_{j,k}\} & =&\{t+j+k,t-1+j+k,|k-j|,1+|k-j| \mid 1 \le j,k \le m \}
\nonumber \\
& \cup &
~~~~~~~~~~~~~~~~~\underbrace{\{0,j,j,j+1,j+1
\mid 1 \le j \le m \}}_{\mbox{omit if $t=0$}}. \nonumber \\
\end{eqnarray}
\begin{thm}
$$H^{\cdot}_{cusp}(SL_n({\Z});{\R}) =0,~1<n<27.$$
\end{thm}
{\bf Proof:} Let $$g_p(x)=\((1-\f{|x|}{p})\cos(\f{\pi
x}{p})+\f{1}{\pi}\sin(\f{\pi|x|}{p})\)/\cosh(x/2).$$  Then
$h_{1,p}(r)=\hat{g}_{1,p}(r)$ is positive in the critical strip
$|\Im r|<\f{1}{2}$. For our cohomological forms $D=1$ at full
level and with the $\mu_{jk}$'s as above we arrive at a
contradiction to the positivity inequality  (see Table
\ref{cohomtable}). \bx

\begin{table}
\begin{tabular}{||||l||l|c|||l||l|c|||l||l|c|||l||l|c||||}
 \hline\hline
\label{thomas}
\em{n} & \em{$t$} & LHS(\ref{polepos}) &
\em{n} & \em{$t$} & LHS(\ref{polepos}) &
\em{n} & \em{$t$} & LHS(\ref{polepos}) &
\em{n} & \em{$t$} & LHS(\ref{polepos}) \\
\hline
2  & 3. & -2.821 & 15 & 6. & -111.4 & 9  & 6. & -71.43 & 22 & 6. &  -77.30 \\
3  & 6. & -8.113 & 16 & 6. & -112.1 & 10 & 6. &  -80.27 & 23 & 6. & -64.06 \\
4  & 6. & -17.02 & 17 & 6. & -112.4 & 11 & 6. & -89.68 & 24 & 6. & -46.70 \\
5  & 6. & -28.30 & 18 & 6. & -109.2 & 12 & 6. &  -96.45 & 25 & 6. & -28.18 \\
6  & 6. & -38.51  & 19 & 6. & -105.4 & 13 & 6. & -103.4 & 26 & 6. & -5.388 \\
7  & 6. & -50.30 & 20 & 6. &  -97.87 & 14 & 6. & -107.5  & & & \\

\hline\hline
\end{tabular}
\caption[Cohomology]{The (numerical) proof of the cohomology theorem.
The left-hand side of (\ref{polepos}),
$\int_{{\R}} g(x)\(e^{x/2}+e^{-x/2}\)dx +
2 \Re \sum_{j=1}^m \sum_{k=1}^m l(\mu_{jk}) + g(0)\log D$
must be positive if the cusp form
exists, and this table shows it is negative for $n<27$.}
\label{cohomtable}
\end{table}

\subsection*{Remarks}
We proved this for $n<23$ in \cite{Miller} with the
Rankin-Selberg
L-functions but without Weil's formula
(instead using the Mittag-Leffler
expansion).  Fermigier \cite{Fermig}
 proved a weaker result using
Weil's formula but with the standard L-function.
The above theorem surpasses both.

Nothing is known about these cuspidal Betti numbers for $n\ge 27$,
let alone if they ever non-zero.

\vspace{1.5 cm}

 {\bf Current Address:}

\vspace{.4 cm}

{\sc  Department of Mathematics

Hill Center-Busch Campus

 Rutgers, The State University of New
Jersey

110 Frelinghuysen Road

Piscataway, NJ 08854-8019

} \vspace{.2 cm}
 {\tt miller@math.rutgers.edu}

\end{document}